\documentclass[a4paper,12pt,twoside,leqno]{article} 
\usepackage[english]{babel} 
\usepackage{amsmath,amssymb,amsthm,amsfonts} 
\usepackage[dvips]{graphicx}
\usepackage[all]{xy}

\usepackage{hyperref}
\hypersetup{
  pdftitle={Steinitz classes of some abelian and nonabelian extensions of even degree},
  pdfauthor={Alessandro Cobbe},
  pdfsubject={Steinitz classes},
  pdfkeywords={Steinitz classes}}

\usepackage{fancyhdr}
\newcommand{\Q}{\mathbb Q}
\newcommand{\N}{\mathbb N}
\newcommand{\Z}{\mathbb Z}

\renewcommand{\P}{\mathfrak P}
\newcommand{\qq}{\mathfrak Q}
\newcommand{\p}{\mathfrak p}
\newcommand{\q}{\mathfrak q}
\newcommand{\gal}{\mathrm{Gal}}
\renewcommand{\epsilon}{\varepsilon}
\newcommand{\disc}{\mathrm{d}}
\newcommand{\norm}{\mathrm{N}}

\newcommand{\st}{\mathrm{st}}

\newcommand{\rt}{\mathrm{R}_t}
\newcommand{\cl}{\mathrm{Cl}}
\newcommand{\aut}{\mathrm{Aut}}

\newcommand{\G}{\mathcal G}
\newcommand{\oo}{\mathcal O}

\newcommand{\comment}[1]{}

\newtheorem{teo}{Theorem}[section]
\newtheorem{lemma}[teo]{Lemma}
\newtheorem{coroll}[teo]{Corollary}
\newtheorem{prop}[teo]{Proposition}
\newtheorem{defn}[teo]{Definition}

\title{Steinitz classes of some abelian and nonabelian extensions of even degree}
\author{Alessandro Cobbe}

\begin{document}
\maketitle

\markboth{Abstract}{Abstract}
\begin{section}*{Abstract}
\addcontentsline{toc}{section}{Abstract}
The Steinitz class of a number field extension $K/k$ is an ideal class in the ring of integers $\oo_k$ of $k$, which, together with the degree $[K:k]$ of the extension determines the $\oo_k$-module structure of $\oo_K$. We call $\rt(k,G)$ the classes which are Steinitz classes of a tamely ramified $G$-extension of $k$. We will say that those classes are realizable for the group $G$; it is conjectured that the set of realizable classes is always a group.

In this paper we will develop some of the ideas contained in \cite{Cobbe} to obtain some results in the case of groups of even order. In particular we show that to study the realizable Steinitz classes for abelian groups, it is enough to consider the case of cyclic groups of $2$-power degree. 
\end{section}

\markboth{Introduction}{Introduction}
\begin{section}*{Introduction}
\addcontentsline{toc}{section}{Introduction}

Let $K/k$ be an extension of number fields and let $\oo_K$ and $\oo_k$ be their rings of integers. By Theorem 1.13 in \cite{Narkiewicz} we know that
\[\oo_K\cong \oo_k^{[K:k]-1}\oplus I\]
where $I$ is an ideal of $\oo_k$. By Theorem 1.14 in \cite{Narkiewicz} the $\oo_k$-module structure of $\oo_K$ is determined by $[K:k]$ and the ideal class of $I$. This class is called the \emph{Steinitz class} of $K/k$ and we will indicate it by $\st(K/k)$. Let $k$ be a number field and $G$ a finite group, then we define:
\[\rt(k,G)=\{x\in\cl(k):\ \exists K/k\text{ tame, }\gal(K/k)\cong G, \st(K/k)=x\}.\]

In this paper we will use the notations and some results from \cite{Cobbe} to study the realizable classes for some groups of even order. In particular we will obtain some explicit results for abelian groups. In this case it follows by \cite{McCulloh_Crelle} that the realizable Steinitz classes form a group, but there is no known explicit characterization of them in the most general situation.

This paper is a slightly shortened version of parts of the author's PhD thesis \cite{tesi}. For earlier results see \cite{Bruche}, \cite{BrucheSodaigui}, \cite{ByottGreitherSodaigui}, \cite{Carter}, \cite{Carter2}, \cite{CarterSodaigui_quaternionigeneralizzati}, \cite{Endo}, \cite{GodinSodaigui_A4}, \cite{GodinSodaigui_ottaedri}, \cite{Long2}, \cite{Long0},  \cite{Massy}, \cite{McCulloh},  \cite{Sodaigui1}, \cite{Sodaigui2} and \cite{Soverchia}.

\end{section}

\markboth{Acknowledgements}{Acknowledgements}
\begin{section}*{Acknowledgements}
\addcontentsline{toc}{section}{Acknowledgements}
I am very grateful to Professor Cornelius Greither and to Professor Roberto Dvornicich for their advice and for the patience they showed, assisting me in the writing of my PhD thesis with a lot of suggestions. I also wish to thank the Scuola Normale Superiore of Pisa, for its role in my mathematical education and for its support during the time I was working on my PhD thesis.
\end{section}

\begin{section}{Preliminary results}
We start recalling the following two fundamental results.
\begin{teo}\label{discriminante}
If $K/k$ is a finite tame Galois extension then 
\[\disc(K/k)=\prod_{\p}\p^{(e_\p-1)\frac{[K:k]}{e_\p}},\]
where $e_\p$ is the ramification index of $\p$.
\end{teo}

\begin{proof}
This follows by Propositions 8 and 14 of chapter III of \cite{Lang}.
\end{proof}

\begin{teo}\label{stdisc}
Assume $K$ is a finite Galois extension of a number field $k$.
\begin{enumerate}
\item[(a)] If its Galois group either has odd order or has a noncyclic $2$-Sylow subgroup then $\disc(K/k)$ is the square of an ideal and this ideal represents the Steinitz class of the extension.
\item[(b)] If its Galois group is of even order with a cyclic $2$-Sylow subgroup and $\alpha$ is any element of $k$ whose square root generates the quadratic subextension of $K/k$ then $\disc(K/k)/\alpha$ is the square of a fractional ideal and this ideal represents the Steinitz class of the extension.
\end{enumerate}
\end{teo}

\begin{proof}
This is a corollary of Theorem I.1.1 in \cite{Endo}. In particular it is shown in \cite{Endo} that in case (b) $K/k$ does have exactly one quadratic subextension.
\end{proof}

Further, considering Steinitz classes in towers of extensions, we will need the following proposition.

\begin{prop}\label{stintermediateextension}
Suppose $K/E$ and $E/k$ are number fields extensions. Then
\[\st(K/k)=\st(E/k)^{[K:E]}\norm_{E/k}(\st(K/E)).\]
\end{prop}

\begin{proof}
This is Proposition I.1.2 in \cite{Endo}.
\end{proof}

We will use a lot of results proved in \cite{Cobbe}. In this section we prove a few facts we will use later.

\begin{lemma}\label{mcdle}
For any $e|m$ the greatest common divisor, for $l|e$, of the integers $(l-1)\frac{m}{e(l)}$ divides $(e-1)\frac{m}{e}$.
\end{lemma}

\begin{proof}
First of all it is clear that we can assume that $m=e$.

Let $I$ be the $\Z$-ideal generated by the $l-1$, for all primes $l|e$. Then $e\equiv 1\pmod I$, since it is the product of prime factors, each one congruent to $1$ modulo $I$.
It follows that for any prime $l\nmid e$, there exists an $l_1|e$, such that the $l$-component of $l_1-1$, which coincides with that of $(l_1-1)\frac{e}{e(l_1)}$, divides that of $e-1$.
Finally, for any $l|e$, $l$ does not divide $(l-1)\frac{e}{e(l)}$.
\end{proof}

\begin{lemma}\label{congruenza}
Let $m,n,x$ be integers. If $x\equiv 1\pmod{m}$ and any prime $q$ dividing $n$ divides also $m$ then
\[x^n\equiv1\pmod{mn}.\]
\end{lemma}

\begin{proof}
Let $n=q_1\dots q_r$ be the prime decomposition of $n$ ($q_i$ and $q_j$ with $i\neq j$ are allowed to be equal). We prove by induction on $r$ that $x^n\equiv 1\pmod{mn}$. If $r=1$, then $mn=mq_1$ must divide $m^{q_1}$ and there exists $b\in\N$ such that
\[x^n=(1+bm)^{q_1}=1+\sum_{i=1}^{q_1-1}\binom{q_1}{i}(bm)^i+(bm)^{q_1}\equiv 1\pmod{mn}.\]
Let us assume that the lemma is true for $r-1$ and prove it for $r$. Since $q_r|m$, as above, for some $c\in \N$ we have
\[\begin{split}x^n&=(1+cmq_1\dots q_{r-1})^{q_r}=1+\sum_{i=1}^{q_r}\binom{q_r}{i}(cmq_1\dots q_{r-1})^i\equiv 1\pmod{mn}.\end{split}\]
\end{proof}

\begin{defn}
Let $K/k$ be a finite abelian extension of number fields. Then we define $W(k,K)$ in the following equivalent ways (the equivalence is shown in \cite{Cobbe}, Proposition 1.10):
\[\begin{split}
&W(k,K)=\{x\in J_k/P_k:\text{$x$ contains infinitely many primes of absolute}\\&\qquad\qquad \qquad \text{degree $1$ splitting completely in $K$}\}\\
&W(k,K)=\{x\in J_k/P_k:\text{$x$ contains a prime splitting completely in $K$}\}\\
&W(k,K)=\norm_{K/k}(J_K)\cdot P_k/P_k.\end{split}\]
\end{defn}
In the case of cyclotomic extensions we will also use the shorter notation $W(k,m)=W(k,k(\zeta_m))$.

\begin{lemma}\label{Wexp}
If $q|n\Rightarrow q|m$ then
$W(k,m)^n\subseteq W(k,mn)$.
\end{lemma}

\begin{proof}
Let $x\in W(k,m)$. By the definition and by Lemma 1.11 of \cite{Cobbe}, $x$ contains a prime ideal $\p$, prime to $mn$ and such that $\norm_{k/\Q}(\p)\in P_\Q^\mathfrak m$, where $\mathfrak m=m\cdot p_\infty$. Then by Lemma \ref{congruenza}, $\norm_{k/\Q}(\p^n)\in P_\Q^\mathfrak n$, with $\mathfrak n=mn\cdot p_\infty$, and it follows from Lemma 1.12 of \cite{Cobbe} that $x^n\in W(k,mn)$.
\end{proof}

\end{section}

\begin{section}{Some general results}
We recall the following definition, from \cite{Cobbe}.
\begin{defn}
We will call a finite group $G$ \emph{good} if the following properties are verified:
\begin{enumerate}
\item For any number field $k$, $\rt(k,G)$ is a group.
\item For any tame $G$-extension $K/k$ of number fields there exists an element $\alpha_{K/k}\in k$ such that:
\begin{enumerate}
\item[(a)] If $G$ is of even order with a cyclic $2$-Sylow subgroup, then a square root of $\alpha_{K/k}$ generates the quadratic subextension of $K/k$; if $G$ either has odd order or has a noncyclic $2$-Sylow subgroup, then $\alpha_{K/k}=1$.
\item[(b)] For any prime $\p$, with ramification index $e_\p$ in $K/k$, the ideal class of
\[\left(\p^{(e_\p-1)\frac{m}{e_\p}-v_\p(\alpha_{K/k})}\right)^\frac{1}{2}\]
is in $\rt(k,G)$. 
\end{enumerate}
\item For any tame $G$-extension $K/k$ of number fields, for any prime ideal $\p$ of $k$ and any rational prime $l$ dividing its ramification index $e_\p$, the class of the ideal
\[\p^{(l-1)\frac{m}{e_\p(l)}}\]
is in $\rt(k,G)$ and, if $2$ divides $(l-1)\frac{m}{e_\p(l)}$, the class of
\[\p^{\frac{l-1}{2}\frac{m}{e_\p(l)}}\]
is in $\rt(k,G)$.
\item $G$ is such that for any number field $k$, for any class $x\in\rt(k,G)$ and any integer $n$, there exists a tame $G$-extension $K$ with Steinitz class $x$ and such that every non trivial subextension of $K/k$ is ramified at some primes which are unramified in $k(\zeta_{n})/k$.
\end{enumerate}
\end{defn}

In \cite{Cobbe} we prove that abelian groups of odd order are good and, more generally, we construct nonabelian good groups by an iteration of direct and semidirect products.

Let $\G$ be a finite group of order $m$, let $H=C(n_1)\times\dots\times C(n_r)$ be an abelian group of order $n$, with generators $\tau_1,\dots,\tau_r$ and with $n_{i+1}|n_i$. Let
\[\mu:\G\to \aut(H)\]
be an action of $\G$ on $H$ and let
\[0\to H\stackrel\varphi\longrightarrow G\stackrel\psi\longrightarrow \G\to 0\]
be an exact sequence of groups such that the induced action of $\G$ on $H$ is $\mu$. We assume that the group $G$ is determined, up to isomorphism, by the above exact sequence and by the action $\mu$. The following well-known proposition shows a class of situations in which our assumption is true.

\begin{prop}[Schur-Zassenhaus, 1937]\label{Gmetaciclico}
If the order of $H$ is prime to the order of $\G$ then $G$ is a semidirect product:
\[G\cong H\rtimes_\mu \G.\]
\end{prop}
\begin{proof}
This is Theorem 7.41 in \cite{Rotman}.
\end{proof}

We are going to study $\rt(k,G)$, considering in particular the case in which the order of $H$ is even.

We also define
\[\eta_G=\begin{cases}1&\text{if $2\nmid n$ or the $2$-Sylow subgroups of $G$ are not cyclic}\\2&\text{if $2|n$ and the $2$-Sylow subgroups of $G$ are cyclic}\end{cases}\]
and in a similar way we define $\eta_H$ and $\eta_\G$. We will always use the letter $l$ only for prime numbers, even if not explicitly indicated.

We say that $(K,k_1,k)$ is of type $\mu$ if $k_1/k$, $K/k_1$ and $K/k$ are Galois extensions with Galois groups isomorphic to $\G$, $H$ and $G$ respectively and such that the action of $\gal(k_1/k)\cong \G$ on $\gal(K/k_1)\cong H$ is given by $\mu$. For any $\G$-extension $k_1$ of $k$ we define $\rt(k_1,k,\mu)$ as the set of those ideal classes of $k_1$ which are Steinitz classes of a tamely ramified extension $K/k_1$ for which $(K,k_1,k)$ is of type $\mu$.

We will repeatedly use the following generalization of the Multiplication Lemma on page 22 in \cite{Endo} by Lawrence P. Endo.

\begin{lemma}\label{multiplicationmetacyclic}
Let $(K_1,k_1,k)$ and $(K_2,k_1,k)$ be extensions of type $\mu$, such that $(\disc(K_1/k_1),\disc(K_2/k_1))=1$ and $K_1/k_1$ and $K_2/k_1$ have no non-trivial unramified subextensions. Then there exists an extension $(K,k_1,k)$ of type $\mu$, such that $K\subseteq K_1K_2$ and for which
\[\st(K/k_1)=\st(K_1/k_1)\st(K_2/k_1).\]
\end{lemma}

\begin{proof}
This is Lemma 2.5 in \cite{Cobbe}.
\end{proof}

Now we recall some further notations introduced in \cite{Cobbe}.

For any $\tau\in H$ we define 
\[\tilde G_{k,\mu,\tau}=\left\{(g_1,g_2)\in \G\times\gal(k(\zeta_{o(\tau)})/k):\ \mu(g_1)(\tau)=\tau^{\nu_{k,\tau}(g_2)}\right\},\]
where $g_2(\zeta_{o(\tau)})=\zeta_{o(\tau)}^{\nu_{k,\tau}(g_2)}$ for any $g_2\in \gal(k(\zeta_{o(\tau)})/k)$,
\[G_{k,\mu,\tau}=\left\{g\in\gal(k(\zeta_{o(\tau)})/k):\ \exists g_1\in\G,\ (g_1,g)\in \tilde G_{k,\mu,\tau}\right\}\]
and $E_{k,\mu,\tau}$ as the fixed field of $G_{k,\mu,\tau}$ in $k(\zeta_{o(\tau)})$.

Given a $\G$-extension $k_1$ of $k$, there is an injection of $\gal(k_1(\zeta_{o(\tau)})/k)$ into $\G\times\gal(k(\zeta_{o(\tau)})/k)$ (defined in the obvious way). We will always identify $\gal(k_1(\zeta_{o(\tau)})/k)$ with its image in $\G\times\gal(k(\zeta_{o(\tau)})/k)$. So we may consider the subgroup
\[\tilde G_{k_1/k,\mu,\tau}=\tilde G_{k,\mu,\tau}\cap \gal(k_1(\zeta_{o(\tau)})/k)\]
of $\tilde G_{k,\mu,\tau}$. 
Let $Z_{k_1/k,\mu,\tau}$ be its fixed field in $k_1(\zeta_{o(\tau)})$.

If $k_1\cap k(\zeta_{o(\tau)})=k$ then $\gal(k_1(\zeta_{o(\tau)})/k)\cong \G\times\gal(k(\zeta_{o(\tau)})/k)$ and hence 
$\tilde G_{k_1/k,\mu,\tau}=\tilde G_{k,\mu,\tau}$.

Now we can state one of the principal theorems proved in \cite{Cobbe}.

\begin{teo}\label{abelcyclicresult}
Let $k$ be a number field and let $\G$ be a good group of order $m$. Let $H=C(n_1)\times\dots\times C(n_r)$ be an abelian group of odd order prime to $m$ and let $\mu$ be an action of $\G$ on $H$. Then
\[\rt(k,H\rtimes_\mu \G)=\rt(k,\G)^n \prod_{l|n}\prod_{\tau\in H(l)\setminus\{1\}}W\left(k,E_{k,\mu,\tau}\right)^{\frac{l-1}{2}\frac{mn}{o(\tau)}}.\]
Furthermore $G=H\rtimes_\mu\G$ is good.
\end{teo}
\begin{proof}
This is Theorem 2.19 of \cite{Cobbe}.
\end{proof}

In this paper we will obtain some results also for abelian groups $H$ of even order, if all the other hypotheses of the above theorem continue to hold.

\begin{lemma}\label{metaconstructivelemma1}
Let $k_1$ be a tame $\G$-extension of $k$ and let $x\in W(k,k_1(\zeta_{n_1}))$. Then there exist tame extensions of $k_1$ of type $\mu$, whose Steinitz classes (over $k_1$) are $\iota(x)^{\eta_H\alpha}$, where
\[\alpha=\sum_{i=1}^{r} \frac{n_i-1}{2}\frac{n}{n_i}+\frac{n_1-1}{2}\frac{n}{n_1}.\]
In particular there exist tame extensions of $k_1$ of type $\mu$ with trivial Steinitz class. 

We can choose these extensions so that they are unramified at all infinite primes, that the discriminants are prime to a given ideal $I$ of $\oo_k$ and that all their proper subextensions are ramified.
\end{lemma}

\begin{proof}
As in the proof of Lemma 2.10 in \cite{Cobbe} we can construct an extension of type $\mu$ with discriminant
\[\disc=\left(\prod_{i=1}^{r}\q_i^{(n_i-1)\frac{n}{n_i}}\right)\q_{r+1}^{(n_1-1)\frac{n}{n_1}}\oo_{k_1},\]
where $\q_1,\dots,\q_{r+1}$ are prime ideals of $k$ in the class of $x$.

If $H$ is of odd order or the $2$-Sylow subgroup of $H$ is cyclic then the result follows immediately by Theorem \ref{stdisc} (a). If this is not the case then by Theorem \ref{stdisc} (b) we obtain extensions whose Steinitz classes have $x^{2\alpha}$ as their square. We may construct infinitely many such $\mu$-extensions  whose discriminants over $k_1$ are relatively prime and so, by the pigeonhole principle, there are two of them, which we call $K_1$ and $K_2$, with the same Steinitz class. Then the conclusion follows considering the extension $K$ given by Lemma \ref{multiplicationmetacyclic}.

As in the proof of Lemma 2.10 in \cite{Cobbe} we can assume that all the additional conditions are also verified.
\end{proof}

\begin{lemma}\label{metaconstructivelemma}
Let $k_1$ be a $\G$-extension of $k$, let $H$ be a group of even order $n$, let $\tau\in H(2)\setminus\{1\}$ and let $x$ be any class in $W(k,Z_{k_1/k,\mu,\tau})$. Then there exist extensions of $k_1$ of type $\mu$, whose Steinitz classes (over $k_1$) are $\iota(x)^{\eta_H\alpha_j}$, where:
\[\alpha_1=\frac{n}{2},\tag{a}\]
\[\alpha_2=(o(\tau)-1)\frac{n}{o(\tau)},\tag{b}\]
Further there exist extensions whose Steinitz classes have $\iota(x)^{2\alpha_j}$ as their square. 
We can choose these extensions so that they satisfy the additional conditions of Lemma \ref{metaconstructivelemma1}.
\end{lemma}

\begin{proof}
\begin{enumerate}
\item[(a)] As in the proof of Lemma 2.11 (a) in \cite{Cobbe} we can construct an extension of type $\mu$ with discriminant
\[\disc(K/k_1)\left((\q_1\q_2)^{\frac{n}{2}}\oo_{k_1}\right),\]
where $\q_1$ and $\q_2$ are prime ideals of $k$ in the class of $x$ and $(K,k_1,k)$ is a $\mu$-extension of $k_1$ with trivial Steinitz class, obtained by Lemma \ref{metaconstructivelemma1}. 
Its Steinitz class has $\iota(x)^{2\alpha_1}$ as its square and we conclude as in Lemma \ref{metaconstructivelemma1}.

\item[(b)] In this case, as in the proof of Lemma 2.11 (b) of \cite{Cobbe}, we obtain an extension of type $\mu$ with discriminant
\[\disc(K/k_1)\left((\q_1\q_2)^{(o(\tau)-1)\frac{n}{o(\tau)}}\oo_{k_1}\right).\]
We conclude as in (a).
\end{enumerate}
\end{proof}

\begin{lemma}\label{mudiretto}
Let $k_1/k$ be a $\G$-extension of number fields, let $H(2)$ be the $2$-Sylow subgroup of $H$ and $\tilde H$ such that $H=H(2)\times \tilde H$. Let $\mu_{\tilde H}$ and $\mu_{H(2)}$ be the actions of $\G$ induced by $\mu$ on $\tilde H$ and $H(2)$ respectively. Then
\[\rt(k_1,k,\mu_{\tilde H})^{n(2)}\subseteq\rt(k_1,k,\mu).\]
\end{lemma}

\begin{proof}
Let $x\in\rt(k_1,k,\mu_{\tilde H})$ and let $(\tilde K,k_1,k)$ be a $\mu_{\tilde H}$-extension of $k_1$ with Steinitz class $x$, which is the class of
\[\disc(\tilde K/k_1)^\frac{1}{2}.\]
Let $(K,k_1,k)$ be a $\mu_{H(2)}$-extension of $k_1$ with trivial Steinitz class and such that $K/k_1$ and $\tilde K/k_1$ are arithmetically disjoint (such an extension exists because of Lemma \ref{metaconstructivelemma1}). The Steinitz class of $K/k$ is the class of
\[\left(\frac{\disc(K/k_1)}{\alpha}\right)^\frac{1}{2}\]
for a certain $\alpha\in k_1$.
Then the extension $(K\tilde K,k_1,k)$ is a $\mu$-extension and its Steinitz class is the class of
\[\left(\frac{\disc(K\tilde K/k)}{\alpha^\frac{n}{n(2)}}\right)^\frac{1}{2}=\disc(\tilde K/k)^\frac{n(2)}{2}\left(\frac{\disc(K/k)}{\alpha}\right)^\frac{n}{2n(2)}\]
which is $x^{n(2)}$.
\end{proof}

At this point we can prove the following proposition.
\begin{prop}\label{metaconstructive}
Let $l\neq 2$ be a prime dividing $n$ and let $\tau\in H(l)$, then
\[\iota\left(W\left(k,Z_{k_1/k,\mu,\tau}\right)\right)^{\frac{l-1}{2}\frac{n}{o(\tau)}}\subseteq \rt(k_1,k,\mu)\]
If $2|n$ then, for any $\tau\in H(2)$,
\[\iota\left(W\left(k,Z_{k_1/k,\mu,\tau}\right)\right)^{\eta_H\frac{n}{o(\tau)}}\subseteq \rt(k_1,k,\mu)\]
and
\[\iota\left(W\left(k,Z_{k_1/k,\mu,\tau}\right)\right)^{2\frac{n}{o(\tau)}}\subseteq \rt(k_1,k,\mu)^2.\]
We can choose the corresponding extensions so that they satisfy the additional conditions of Lemma \ref{metaconstructivelemma1}.
\end{prop}

\begin{proof}
The first inclusion follows immediately by Lemma 2.12 of \cite{Cobbe} and by Lemma \ref{mudiretto}.

Now let us assume that $2|n$, let $\tau\in H(2)$ and let $x\in W(k,Z_{k_1/k,\mu,\tau})$. It follows from Lemma \ref{multiplicationmetacyclic} and Lemma \ref{metaconstructivelemma} that $\iota(x)^{\eta_H\beta_2}$ is in $\rt(k_1,k,\mu)$ and $\iota(x)^{2\beta_2}$ is in $\rt(k_1,k,\mu)^2$, where
\[\beta_2=\gcd\left(\frac{n}{2},(o(\tau)-1)\frac{n}{o(\tau)}\right).\]
So we obtain
\[\iota(x)^{\eta_H\frac{n}{o(\tau)}}\in\rt(k_1,k,\mu)\]
and
\[\iota(x)^{2\frac{n}{o(\tau)}}\in\rt(k_1,k,\mu)^2.\]
To conclude we observe that applying Lemma \ref{multiplicationmetacyclic} to extensions $(K_1,k_1,k)$ and $(K_2,k_1,k)$ of type $\mu$ which satisfy the additional conditions of Lemma \ref{metaconstructivelemma1}, we obtain an extension $(K,k_1,k)$ which still satisfies the same conditions.
\end{proof}

\begin{prop}\label{abelcyclicresultoneinclusion}
Let $a$ be a multiple of a positive integer $n_1$. Let $k$ be a number field and let $\G$ be a finite group such that for any class $x\in\rt(k,\G)$ there exists a tame $\G$-extension $k_1$ with Steinitz class $x$ and such that every subextension of $k_1/k$ is ramified at some primes which are unramified in $k(\zeta_{a})/k$.

Let $H=C(n_1)\times\dots\times C(n_r)$ be an abelian group of order $n$ and let $\mu$ be an action of $\G$ on $H$. We assume that the exact sequence
\[0\to H\stackrel\varphi\longrightarrow G\stackrel\psi\longrightarrow \G\to 0,\]
in which the induced action of $\G$ on $H$ is $\mu$, determines the group $G$, up to isomorphism. Further we assume that $H$ is of odd order or with noncyclic $2$-Sylow subgroup, or that $\G$ is of odd order. Then
\[\rt(k,H\rtimes_\mu \G)\supseteq\rt(k,\G)^n \prod_{\substack{l|n\\l\neq 2}}\prod_{\tau\in H(l)}\!W\left(k,E_{k,\mu,\tau}\right)^{\frac{l-1}{2}\frac{mn}{o(\tau)}}\!\!\prod_{\tau\in H(2)}\!W\left(k,E_{k,\mu,\tau}\right)^{\frac{\eta_Gmn}{o(\tau)}}.\]

Further we can choose tame $G$-extensions $K/k$ with a given Steinitz class (of the ones considered above), such that every nontrivial subextension of $K/k$ is ramified at some primes which are unramified in $k(\zeta_a)/k$.
\end{prop}
\begin{proof}
Let $x\in\rt(k,\G)$ and let $k_1$ be a tame $\G$-extension of $k$, with Steinitz class $x$, and such that every subextension of $k_1/k$ is ramified at some primes which are unramified in $k(\zeta_{a})/k$. Thus, since $a$ is a multiple of $n_1$, it follows also that $k_1\cap k(\zeta_{n_1})=k$. 

By Lemma \ref{multiplicationmetacyclic}, Lemma 2.4 in \cite{Cobbe}, Proposition \ref{metaconstructive} and Proposition \ref{stintermediateextension} we obtain
\[\rt(k,H\rtimes_\mu \G)\supseteq x^n \prod_{\substack{l|n\\l\neq 2}}\prod_{\tau\in H(l)}W\left(k,E_{k,\mu,\tau}\right)^{\frac{l-1}{2}\frac{mn}{o(\tau)}}\prod_{\tau\in H(2)}W\left(k,E_{k,\mu,\tau}\right)^{\frac{\eta_H mn}{o(\tau)}},\]
from which we obtain the result we wanted to prove, if $\eta_H=\eta_G$. 

With our hypotheses $\eta_H\neq\eta_G$ implies that the order of $H$ is odd, i.e. that there does not exist any nontrivial $\tau\in H(2)$. Hence we obtain the desired result also in this case.
\end{proof}

\begin{prop}\label{semidiretto2}
Let $\tau,\tilde\tau\in H(2)$ be elements such that $\tau,\tilde\tau,\tau\tilde\tau$ are all of the same order. Let $k_1$ be a $\G$-extension of $k$. Then
\[\iota(W(k,Z_{k_1/k,\mu,\tau}Z_{k_1/k,\mu,\tilde\tau}Z_{k_1/k,\mu,\tau\tilde\tau}))^\frac{n}{2o(\tau)}\subseteq \rt(k_1,k,\mu).\]
In particular, if $Z_{k_1/k,\mu,\tau}=Z_{k_1/k,\mu,\tau}Z_{k_1/k,\mu,\tilde\tau}Z_{k_1/k,\mu,\tau\tilde\tau}$, the factor\footnote{If the order of $\tau$ is $2$ or $4$ this condition is obviously verified (eventually, after renaming $\tau$, $\tilde\tau$ and $\tau\tilde\tau$).}
\[W(k,E_{k,\mu,\tau})^\frac{mn}{2o(\tau)}\]
can be added in the right hand side of the expression of Proposition \ref{abelcyclicresultoneinclusion}, giving more realizable classes. The additional condition of Proposition \ref{abelcyclicresultoneinclusion} is also satisfied.
\end{prop}

\begin{proof}
Let 
\[x\in W(k,Z_{k_1/k,\mu,\tau}Z_{k_1/k,\mu,\tilde\tau}Z_{k_1/k,\mu,\tau\tilde\tau}).\]

We will use all the notations of the proof of Lemma 2.11 of \cite{Cobbe} and we also consider prime ideals $\q_1,\q_2,\q_3$ with analogous conditions.

We define $\varphi_i:\kappa_{\qq_{i}}^*\to H$, for $i=1,2,3$, posing
\[\varphi_1(g_{\qq_1})=\tau,\]
\[\varphi_2(g_{\qq_2})=\tilde\tau,\]
and
\[\varphi_3(g_{\qq_3})=(\tau\tilde\tau)^{-1}.\]
In the usual way we obtain an extension of type $\mu$ with discriminant
\[\disc(K/k_1)\left((\q_1\q_2\q_3)^{(o(\tau)-1)\frac{n}{o(\tau)}}\oo_{k_1}\right)\]
and Steinitz class $\iota(x)^{\alpha_3}$ (with the above hypotheses the $2$-Sylow subgroup of $H$ can not be cyclic), where
\[\alpha_3=3(o(\tau)-1)\frac{n}{2o(\tau)}.\]
Thus by Lemma \ref{multiplicationmetacyclic} and Proposition \ref{metaconstructive} we obtain that
\[\iota(W(k,Z_{k_1/k,\mu,\tau}Z_{k_1/k,\mu,\tilde\tau}Z_{k_1/k,\mu,\tau\tilde\tau}))^\frac{n}{2o(\tau)}\subseteq \rt(k_1,k,\mu).\]
To prove that \[W(k,E_{k,\mu,\tau})^\frac{mn}{2o(\tau)}\]
can be added in the expression of Proposition \ref{abelcyclicresultoneinclusion}, it is now enough to use Lemma 2.4 in \cite{Cobbe}, assuming that $k_1\cap k(\zeta_{o(\tau)})=k$ and that every subextension of $k_1/k$ is ramified (we can make these assumptions thanks to the hypotheses of Proposition \ref{metaconstructive}).
\end{proof}

\begin{lemma}\label{primoramificatosopra}
Let $(K,k_1,k)$ be a tame $\mu$-extension and let $\P$ be a prime in $k_1$ ramifying in $K/k_1$ and let $\p$ be the corresponding prime in $k$. Then
\[x\in W(k,Z_{k_1/k,\mu,\tau})\subseteq W(k,E_{k,\mu,\tau})\subseteq \bigcap_{l|e_\P}W(k,E_{k,\mu,\tau(l)})\]
where $x$ is the class of $\p$ and $\tau$ generates $([U_\P],K/k_1)$.
\end{lemma}
\begin{proof}
This is Lemma 2.14 in \cite{Cobbe}\footnote{In this case we made no assumption concerning the parity of the order of $H$, so the result holds also in the present setting.}.
\end{proof}

\begin{lemma}\label{metaciclicoanalitico2}
Let $\G$ be a good group, let $H$ be an abelian group of order prime to that of $\G$, with trivial or noncyclic $2$-Sylow subgroup, and let $\mu$ be an action of $\G$ on $H$.
Suppose $(K,k_1,k)$ is tamely ramified and of type $\mu$.
Let $e_\p$ be the ramification index of a prime $\p$ in $k_1/k$ and $e_\P$ be the ramification index of a prime $\P$ of $k_1$ dividing $\p$ in $K/k_1$. Then the class of
\[\left(\p^{(e_\p e_\P-1)\frac{mn}{e_\p e_\P}-v_\p(\alpha_{k_1/k}^n)}\right)^\frac{1}{2}\]
is in
\[ \rt(k,\G)^n\cdot\prod_{l|n}\prod_{\tau\in H(l)\setminus\{1\}}W\left(k,E_{k,\mu,\tau}\right)^{\frac{l-1}{2}\frac{mn}{o(\tau)}}.\]
\end{lemma}

\begin{proof}
If the order of $H$ is odd, then this is Lemma 2.17 in \cite{Cobbe}. So we can assume that the order of $H$ is even and the order of $\G$ is odd. Exactly as in the proof of Lemma 2.17 in \cite{Cobbe} we obtain
\[\p^{(e_\p e_\P-1)\frac{mn}{e_\p e_\P}}=\p^{a_\p(e_\p-1)\frac{mn}{e_\p}}\prod_{l|e_\P} \p^{b_{\p,l}(l-1)\frac{mn}{e_\P(l)}}. 
\]

By our assumption $\alpha_{k_1/k}=1$ and we conclude by the hypothesis that $\G$ is good, by Lemma \ref{primoramificatosopra} and by the fact that, for any prime $l$ dividing $e_\P$, $(l-1)\frac{mn}{e_\P(l)}$ is even (in the case $l=2$ this is due to the fact that the inertia group at $\P$ must be cyclic, while the $2$-Sylow subgroup of $H$ is not).
\end{proof}

\begin{lemma}\label{metaciclicoanalitico1}
Under the same hypotheses as in the preceding lemma, if $l|e_\p e_\P$, the class of
\[\p^{(l-1)\frac{mn}{e_\p(l)e_\P(l)}}\]
is in
\[\rt(k,\G)^{n}\prod_{\tau\in H(l)\setminus\{1\}}W\left(k,E_{k,\mu,\tau}\right)^{\frac{l-1}{2}\frac{mn}{o(\tau)}}.\]
and, if $2$ divides $(l-1)\frac{mn}{e_\p(l)e_\P(l)}$, the class of
\[\p^{\frac{l-1}{2}\frac{mn}{e_\p(l)e_\P(l)}}\]
is in
\[\rt(k,\G)^{n}\prod_{\tau\in H(l)\setminus\{1\}}W\left(k,E_{k,\mu,\tau}\right)^{\frac{l-1}{2}\frac{mn}{o(\tau)}}.\]
\end{lemma}

\begin{proof}
As in the previous lemma we can assume that the order of $H$ is even, since the odd case has been proved in Lemma 2.18 of \cite{Cobbe}. Thus if $l$ is a prime dividing $e_\p$, then $l$ is odd and so $2$ divides $(l-1)\frac{m}{e_\p(l)}$ and the class of
\[\p^{\frac{l-1}{2}\frac{m}{e_\p(l)}}\]
is in $\rt(k,\G)$, by the hypothesis that $\G$ is good. We conclude that the class of
\[\p^{\frac{l-1}{2}\frac{mn}{e_\p(l)e_\P(l)}}=\p^{\frac{l-1}{2}\frac{mn}{e_\p(l)}}\]
is in $\rt(k,\G)^n$. 

If $l$ divides $e_\P$, then $(l-1)\frac{mn}{e_\p(l)e_\P(l)}$ is even (by hypothesis the $2$-Sylow subgroup of $H$ is not cyclic and thus $\frac{n}{e_\P(2)}$ is even). We conclude by Lemma \ref{primoramificatosopra} that the class of
\[\p^{\frac{l-1}{2}\frac{mn}{e_\p(l)e_\P(l)}}=\p^{\frac{l-1}{2}\frac{mn}{e_\P(l)}}\]
is in $W\left(k,E_{k,\mu,\tau}\right)^{\frac{l-1}{2}\frac{mn}{o(\tau)}}$ for some $\tau\in H(l)\setminus\{1\}$.
\end{proof}

\begin{prop} \label{2vectorgood}
Let $k$ be a number field and let $\G$ be a good group of odd order.

Let $n>1$ be an integer, let $H=C(2)^{(n)}=C(2)\times\dots\times C(2)$ and let $\mu$ be an action of $\G$ on $H$. Then
\[\rt(k,H\rtimes_\mu \G)=\rt(k,\G)^{2^n} \cl(k)^{m2^{n-2}}.\]
Further $G=H\rtimes_\mu\G$ is good.
\end{prop}

\begin{proof}
Clearly $E_{k,\mu,\tau}=k$, i.e. $W(k,E_{k,\mu,\tau})=\cl(k)$ for any $\tau\in H(2)=H$. Thus, by Propositions \ref{abelcyclicresultoneinclusion} and \ref{semidiretto2},
\[\rt(k,H\rtimes_\mu \G)\supseteq\rt(k,\G)^{2^n} \cl(k)^{m2^{n-2}}.\]
The opposite inclusion comes from Theorems \ref{discriminante} and \ref{stdisc} and from Lemma \ref{metaciclicoanalitico2}. So we obtain an equality and, in particular, this gives the first property of good groups. The other properties follow now respectively from Lemmas \ref{metaciclicoanalitico2} and \ref{metaciclicoanalitico1} and from Propositions \ref{abelcyclicresultoneinclusion} and \ref{semidiretto2}.
\end{proof}

If $\G$ is cyclic of order $2^n-1$ and the representation $\mu$ is faithful, then the above proposition is one of the results proved by Nigel P. Byott, Cornelius Greither and Boucha\"ib Soda\"igui in \cite{ByottGreitherSodaigui}.

\vspace{10pt}

\textbf{Example.} The group $A_4$, which is isomorphic to a semidirect product of the form $(C(2)\times C(2))\rtimes_\mu C(3)$, is good by Proposition \ref{2vectorgood}. We calculate its realizable classes:
\[\rt(k,A_4)=W(k,3)^4\cl(k)^3\supseteq\cl(k)^8\cl(k)^3=\cl(k)\]
and hence
\[\rt(k,A_4)=\cl(k).\]
This result has been obtained by Marjory Godin and Boucha\"ib Soda\"igui in \cite{GodinSodaigui_A4}.

\end{section}

\begin{section}{Abelian extensions of even degree}
We will conclude this paper considering the case of abelian groups of even order. To this aim we will use the preceding results and notations, with the assumption that $\G$ is the trivial group, i.e. that $G=H$.

It follows by the paper \cite{McCulloh_Crelle} of Leon McCulloh that $\rt(k,G)$ is a group for any finite abelian group $G$. However, this result does not yield an explicit description of $\rt(k,G)$, which is known only if the order of $G$ is odd and in a few other cases.

\begin{teo}\label{oddabelian}
Let $k$ be a number field and let $G=C(n_1)\times\dots\times C(n_r)$ with $n_{i+1}|n_i$ be an abelian group of odd order. Then
\[\rt(k,G)=\prod_{l|n} W(k,n_1(l))^{\frac{l-1}{2}\frac{n}{n_1(l)}}.\]
\end{teo}

\begin{proof}
This result was proved by Endo in his PhD thesis \cite{Endo} (in a slightly different form), but it is also a particular case of Theorem 2.19 of \cite{Cobbe}, using also Lemma \ref{Wexp}.
\end{proof}

Further we will also use the following proposition proved by Endo.
\begin{prop}\label{W2}
For any number field $k$
\[W(k,2^n)\subseteq \rt(k,C(2^n)).\]
\end{prop}
\begin{proof}
This is Proposition II.2.4 in \cite{Endo}.
\end{proof}

The equality of Theorem \ref{oddabelian} is not true in general for abelian groups of even order. Nevertheless it is not difficult to prove one inclusion.

\begin{prop}\label{analytisch}
Let $k$ be a number field and let $G$ be an abelian group, then
\[\rt(k,G)\subseteq\prod_{l|n} W(k,n_1(l))^{\frac{l-1}{2}\frac{n}{n_1(l)}}.\]
\end{prop}

\begin{proof}
Let $K/k$ be a tamely ramified extension of number fields with Galois group $G$.
By Theorem \ref{discriminante} and by Lemma \ref{mcdle} there exist $b_{e_\p,l}\in\Z$ such that
\[\begin{split}\disc(K/k)&=\prod_{e_\p\neq 1}\p^{(e_\p-1)\frac{n}{e_\p}}=\prod_{e_\p\neq 1}\prod_{l|e_\p}\p^{b_{e_\p,l}(l-1)\frac{n}{e_\p(l)}}=\prod_{l|n}\prod_{e_\p(l)\neq 1}\p^{b_{e_\p,l}(l-1)\frac{n}{e_\p(l)}}.\end{split}\]
Since $K/k$ is tame, the ramification index $e_\p$ of a prime $\p$ in $K/k$ divides $n_1$. Thus, defining
\[J_{l}=\prod_{e_\p(l)\neq 1}\p^{b_{e_\p,l}\frac{n_1(l)}{e_\p(l)}},\]
we obtain
\[\disc(K/k)=\prod_{l|n}J_l^{(l-1)\frac{n}{n_1(l)}}\]
and by Lemma 1.13 of \cite{Cobbe} and Lemma \ref{Wexp} the class of the ideal $J_l$ belongs to $W(k,n_1(l))$.
We easily conclude by Theorem \ref{stdisc}.
\end{proof}

\begin{prop}\label{konstruktiv}
Let $l\neq 2$ be a prime dividing $n$, then
\[W(k,n_1(l))^{\frac{l-1}{2}\frac{n}{n_1(l)}}\subseteq \rt(k,G).\]
If $2|n$ then
\[W(k,n_1(2))^{\eta_G\frac{n}{n_1(2)}}\subseteq \rt(k,G)\]
and
\[W(k,n_1(2))^{2\frac{n}{n_1(2)}}\subseteq \rt(k,G)^2.\]
We can choose the corresponding extensions so that they satisfy the additional conditions of Lemma \ref{metaconstructivelemma1}.
\end{prop}

\begin{proof}
This is a particular case of Proposition \ref{metaconstructive}.
\end{proof}

\begin{lemma}\label{pezziuguali}
If $2|n$ and $n_2(2)\neq 1$ then
\[W(k,n_2(2))^{\frac{n}{2n_2(2)}}\subseteq\rt(k,G).\]
We can choose the corresponding extensions so that they satisfy the additional conditions of Lemma \ref{metaconstructivelemma1}.
\end{lemma}

\begin{proof}
This is a particular case of Proposition \ref{semidiretto2}.
\end{proof}

Using this lemma we can easily prove a first interesting proposition, which gives a characterization of realizable classes in a particular situation.

\begin{prop}\label{abelianiparibelli}
Let $k$ be a number field, let $G=C(n_1)\times\dots\times C(n_r)$, with $n_{i+1}|n_i$, be an abelian group of order $n$. If $2|n$ and $n_1(2)=n_2(2)$, then
\[\rt(k,G)=\prod_{l|n} W(k,n_1(l))^{\frac{l-1}{2}\frac{n}{n_1(l)}}\]
and the group $G$ is good.
\end{prop}

\begin{proof}
One inclusion is Proposition \ref{analytisch}.

The other inclusion follows by Proposition \ref{konstruktiv} and Lemma \ref{pezziuguali}, using Lemma \ref{multiplicationmetacyclic}.

Thus, in particular, the first and the fourth property of good groups are verified. Now let $K/k$ be a tamely ramified extension of number fields with Galois group $G$. By Lemma \ref{mcdle} there exist $b_{e_\p,l}\in\Z$ such that
\[\begin{split}\p^{(e_\p-1)\frac{m}{e_\p}}=\prod_{l|e_\p}\p^{b_{e_\p,l}(l-1)\frac{m}{e_\p(l)}}=\prod_{l|e_\p}\p^{\frac{m_1(l)}{e_\p(l)}b_{e_\p,l}(l-1)\frac{m}{m_1(l)}}.\end{split}\]
By Lemma 1.13 of \cite{Cobbe} and Lemma \ref{Wexp}, the class of the ideal $\p^\frac{m_1(l)}{e_\p(l)}$ is contained in $W(k,m_1(l))$. Since $(l-1)\frac{m}{m_1(l)}$ is even for any prime $l$ dividing $e_\p$, we easily conclude that also the second and the third property of good groups hold for $G$.
\end{proof}

\begin{coroll}
Under the assumption of Proposition \ref{abelianiparibelli}, $\rt(k,G)$ is a subgroup of the ideal class group of the number field $k$.
\end{coroll}

The above corollary follows also from \cite{McCulloh_Crelle}, which is however much less explicit than Proposition \ref{abelianiparibelli}. The above description of $\rt(k,G)$ generalizes the result concerning the group $G=C(2)\times C(2)$ in \cite{Sodaigui1}.

\begin{coroll}\label{D2nbuono}
If $n$ is odd then $D_{2n}$ is a good group, it is isomorphic to a semidirect product of the form
\[C(n)\rtimes_\mu(C(2)\times C(2))\]
and
\[\rt(k,D_{2n})=\cl(k)^n\cdot\prod_{l|n}\prod_{\tau\in H(l)\setminus\{1\}}W\left(k,E_{k,\mu,\tau}\right)^{(l-1)\frac{2n}{o(\tau)}}.\]
\end{coroll}
\begin{proof}
It is easy to see that
\[D_{2n}\cong D_n\times C(2)\cong C(n)\rtimes_\mu (C(2)\times C(2)),\]
for a certain action $\mu:C(2)\times C(2)\to \aut(C(n))$. By the above proposition $C(2)\times C(2)$ is good and
\[\rt(k,C(2)\times C(2))=\cl(k).\]
Thus we conclude by Proposition \ref{abelcyclicresult} that $D_{2n}$ is good and we obtain the desired expression for $\rt(k,D_{2n})$.
\end{proof}

An analogous result for a dihedral group $D_n$, where $n$ is an odd integer, is given in Theorem 2.26 of \cite{Cobbe}.

\begin{lemma}\label{dalleciclichealleabeliane}
Let $k$ be a number field and $G=C(n_1)\times\dots\times C(n_r)$ an abelian group of even order $n$. Then
\[\rt(k,C(n_1(2)))^{\frac{n}{n_1(2)}}\subseteq \rt(k,G),\]
where $n_1(2)$ is the maximal power of $2$ dividing $n_1$.
\end{lemma}
\begin{proof}
By hypothesis $G=C(n_1(2))\times \tilde G$, where $\tilde G$ is an abelian group. Let $x\in \rt(k,C(n_1(2)))$ and let $L$ be a tame $C(n_1(2))$-extension whose Steinitz class is $x$. Because of Lemma \ref{metaconstructivelemma1} there exists a tame $\tilde G$-extension $K$ of $k$ whose discriminant is prime to that of $L$ over $k$, with trivial Steinitz class and with no unramified subextensions. The composition of the two extensions is a $G$-extension and its discriminant is
\[\disc(L/k)^\frac{n}{n_1(2)}\disc(K/k)^{n_1(2)}.\]
If the $2$-Sylow subgroup of $G$ is not cyclic then the Steinitz class is the class of
\[\disc(KL/k)^\frac{1}{2}=\disc(L/k)^\frac{n}{2n_1(2)}\disc(K/k)^{n_1(2)/2},\]
that is
\[(x^2)^\frac{n}{2n_1(2)}=x^\frac{n}{n_1(2)}.\]
Now we have to consider the case in which the $2$-Sylow subgroup of $G$ is cyclic. The subextension $k(\sqrt\alpha)$ of $L$ of degree $2$ over $k$ is also a subextension of $KL$. We have $k(\sqrt\alpha)=k\left(\sqrt{\alpha^\frac{n}{n_1(2)}}\right)$ (the exponent $\frac{n}{n_1(2)}$ is odd) and so the Steinitz class of $KL/k$ is the class of the square root of
\[\frac{\disc(KL/k)}{\alpha^\frac{n}{n_1(2)}}=\left(\frac{\disc(L/k)}{\alpha}\right)^\frac{n}{n_1(2)}\disc(K/k)^{n_1(2)},\]
that is exactly $x^\frac{n}{n_1(2)}$.
\end{proof}

\begin{lemma}\label{legameciclichenon}
Let $k$ be a number field and $G=C(n_1)\times\dots\times C(n_r)$ an abelian group of even order $n$, with $n_{i+1}|n_i$ and $n_2(2)\neq 1$. Then
\[\rt(k,G)\subseteq \rt(k,C(n_1(2)))^{\frac{n}{n_1(2)}}\cdot W\left(k,n_2(2)\right)^{\frac{n}{2n_2(2)}}\cdot\prod_{\substack{l|n\\l\neq 2}}W(k,n_1(l))^{\frac{l-1}{2}\frac{n}{n_1(l)}}.\]
\end{lemma}

\begin{proof}
Let $K/k$ be a $G$-Galois extension whose Steinitz class is $x\in\rt(k,G)$ and let $L$ be a subextension of $K/k$ whose Galois group over $k$ is the first component of the $2$-Sylow subgroup $C(n_1(2))\times\dots\times C(n_r(2))$ of $G$. By Theorem 1.3 of \cite{Cobbe} and Proposition II.3.3 in \cite{Neukirch}
\[\begin{split}&e_{\p,K}=e_{\p,K}(2)e_{\p,K}'=\#([U_\p],K/k);\\
&e_{\p,L}=e_{\p,L}(2)=\#([U_\p],L/k)=\#([U_\p],K/k)|_L,\end{split}\]
where $e_{\p,L}$ and $e_{\p,K}$ are the ramification indices of $\p$ in $L$ and $K$ respectively and $e_{\p,K}'$ is odd. By Theorem \ref{discriminante} and Theorem \ref{stdisc}, $x$ is the class of
\[\prod_\p \p^{\frac{e_{\p,K}-1}{2}\frac{n}{e_{\p,K}}}.\]

The class $x_1$ of the ideal
\[\prod_\p \p^{\frac{e_{\p,L}-1}{2}\frac{n}{e_{\p,L}}}\]
is the $n/n_1(2)$-th power of the Steinitz class of $L/k$ and thus
\[x_1\in \rt(k,C(n_1(2)))^\frac{n}{n_1(2)}.\]

Since $e_{\p,L}|e_{\p,K}(2)$ and $2e_{\p,K}(2)|n$ we can define $x_2$ as the class of \[\begin{split}\prod_\p \p^{\left(\frac{e_{\p,K}(2)}{e_{\p,L}}-1\right)\frac{n}{2e_{\p,K}(2)}}=\prod_\p \p^{\left(\frac{e_{\p,K}(2)}{e_{\p,L}}-1\right)\frac{n_2(2)}{e_{\p,K}(2)}\frac{n}{2n_2(2)}}.\end{split}\]
The only primes for which we obtain a nontrivial contribution are those for which $e_{\p,L}<e_{\p,K}(2)$ and for those we must have $e_{\p,K}(2)|n_2(2)$ (since $e_{\p,K}(2)$ must then be the order of a cyclic subgroup of $C(n_2(2))\times\dots\times C(n_{r}(2))$) and thus, recalling Lemma 1.13 of \cite{Cobbe} and Lemma \ref{Wexp},
\[x_2\in W(k,n_2(2))^{\frac{n}{2n_2(2)}}.\]

Let $x_3$ be the class of
\[\prod_\p \p^{\frac{e_{\p,K}'-1}{2}\frac{n}{e_{\p,K}}}=\prod_\p \p^{a_\p\frac{e_{\p,K}'-1}{2}\frac{n}{e_{\p,K}'}}\prod_\p \p^{b_\p\frac{e_{\p,K}'-1}{2}\frac{n}{e_{\p,K}(2)}},\]
where $a_\p$ and $b_\p$ are integers such that
\[\frac{n}{e_{\p,K}}= a_\p\frac{n}{e_{\p,K}'}+b_\p\frac{n}{e_{\p,K}(2)}.\]

By Lemma \ref{mcdle} there exist $b_{\p,l}\in\Z$ such that
\[\prod_\p \p^{a_\p\frac{e_{\p,K}'-1}{2}\frac{n}{e_{\p,K}'}}=\prod_{\substack{l|n\\l\neq 2}}\prod_\p \p^{b_{\p,l}\frac{n_1(l)}{e_{\p,K}'(l)}\frac{l-1}{2}\frac{n}{n_1(l)}}\]
and thus by Lemma 1.13 of \cite{Cobbe} and Lemma \ref{Wexp} the class of this ideal is in
\[\prod_{\substack{l|n\\l\neq 2}}W(k,n_1(l))^{\frac{l-1}{2}\frac{n}{n_1(l)}}.\]
By the same lemmas the class of
\[\prod_\p \p^{b_\p\frac{e_{\p,K}'-1}{2}\frac{n}{e_{\p,K}(2)}}\]
is in
\[W(k,n_1(2))^\frac{n}{n_1(2)},\]
which is contained in
\[\rt(k,C(n_1(2)))^\frac{n}{n_1(2)}\]
by Proposition \ref{W2}. Hence
\[x_3\in \prod_{\substack{l|n\\l\neq 2}}W(k,n_1(l))^{\frac{l-1}{2}\frac{n}{n_1(l)}}\rt(k,C(n_1(2)))^\frac{n}{n_1(2)}.\]

By an easy calculation
\[\frac{e_{\p,K}-1}{2}\frac{n}{e_{\p,K}}=\frac{e_{\p,L}-1}{2}\frac{n}{e_{\p,L}}+\left(\frac{e_{\p,K}(2)}{e_{\p,L}}-1\right)\frac{n}{2e_{\p,K}(2)}+\frac{e_{\p,K}'-1}{2}\frac{n}{e_{\p,K}}\]
and we conclude that $x=x_1x_2x_3$, obtaining the desired inclusion.
\end{proof}

\begin{teo}\label{abelpari2nc}
Let $k$ be a number field, let $G=C(n_1)\times\dots\times C(n_r)$, with $n_{i+1}|n_i$, be an abelian group of order $n$. If $2|n$ and $n_2(2)\neq1$ then
\[\rt(k,G)= \rt(k,C(n_1(2)))^{\frac{n}{n_1(2)}}\cdot W(k,n_1(2))^{\frac{n}{2n_2(2)}}\cdot\prod_{\substack{l|n\\l\neq 2}}W(k,n_1(l))^{\frac{l-1}{2}\frac{n}{n_1(l)}}.\]
\end{teo}

\begin{proof}
\begin{enumerate}
\item[$\subseteq$] This is Lemma \ref{legameciclichenon}.
\item[$\supseteq$] This follows by Proposition \ref{konstruktiv}, by Lemma \ref{pezziuguali} and by Lemma \ref{dalleciclichealleabeliane}, using Lemma \ref{multiplicationmetacyclic}.
\end{enumerate}
\end{proof}

\textbf{Remark. } The only unknown term in the expression for $\rt(k,G)$ in the above theorem is $\rt(k,C(n_1(2)))$. But we really need to determine only its square, because it appears with an even exponent. This simplifies the problem, because this allows us to consider directly the discriminants of the extensions.

\vspace{10pt}

In the second part of the section we consider the case in which the $2$-Sylow subgroup of $G$ is cyclic, i.e. $2|n$ and $n_2(2)=1$.

\begin{lemma}\label{legameciclichenonr11}
If the $2$-Sylow subgroup of $G$ is cyclic, i.e. $2|n$ and $n_2(2)=1$, then
\[\rt(k,G)\subseteq \rt(k,C(n_1(2)))^{\frac{n}{n_1(2)}}\cdot \prod_{\substack{l|n\\l\neq 2}}W(k,n_1(l))^{\frac{l-1}{2}\frac{n}{n_1(l)}}.\]
\end{lemma}

\begin{proof}
Let $K/k$ be a $G$-Galois extension whose Steinitz class is $x\in\rt(k,G)$ and let $L$ be the subextension of $K/k$ whose Galois group over $k$ is the $2$-Sylow subgroup $C(n_1(2))$ of $G$. By Theorem 1.3 of \cite{Cobbe} and Proposition II.3.3 in \cite{Neukirch}
\[\begin{split}&e_{\p,K}=e_{\p,K}(2)e_{\p,K}'=\#([U_\p],K/k);\\
&e_{\p,L}=e_{\p,L}(2)=\#([U_\p],L/k)=\#([U_\p],K/k)|_L,\end{split}\]
where $e_{\p,L}$ and $e_{\p,K}$ are the ramification indices of $\p$ in $L$ and $K$ respectively, $e_{\p,K}'$ is odd and $e_{\p,K}(2)=e_{\p,L}(2)$. Let $\alpha\in k$ be such that $k\subsetneq k(\sqrt\alpha)\subseteq L$.

Since $k(\sqrt\alpha)=k\left(\sqrt{\alpha^{n/n_1(2)}}\right)$, by Theorem \ref{discriminante} and Theorem \ref{stdisc}, $x$ is the class of
\[\left(\frac{\prod_\p \p^{(e_{\p,K}-1)\frac{n}{e_{\p,K}}}}{\alpha^\frac{n}{n_1(2)}}\right)^\frac12.\]

As in the proof of Lemma \ref{legameciclichenon} we can define\footnote{The analogous element in Lemma \ref{legameciclichenon} was called $x_3$.}
\[x_1\in\rt(k,C(n_1(2)))^{\frac{n}{n_1(2)}}\cdot \prod_{\substack{l|n\\l\neq 2}}W(k,n_1(l))^{\frac{l-1}{2}\frac{n}{n_1(l)}}.\]
as the class of the ideal
\[\prod_{\p}\p^{\frac{e_{\p,K}'-1}{2}\frac{n}{e_{\p,K}}}.\]

By Theorem \ref{discriminante} and Theorem \ref{stdisc}, 
\[\left(\frac{\prod_\p \p^{(e_{\p,L}-1)\frac{n_1(2)}{e_{\p,L}}}}{\alpha}\right)^\frac{n}{2n_1(2)}\]
is an ideal, whose class $x_2$ is the $n/n_1(2)$-th power of the Steinitz class of $L/k$. Thus
\[x_2\in\rt(k,C(n_1(2)))^\frac{n}{n_1(2)}.\]
By an easy calculation
\[\left(\frac{\prod_{\p}\p^{(e_{\p,K}-1)\frac{n}{e_{\p,K}}}}{\alpha^\frac{n}{n_1(2)}}\right)^\frac12=\prod_{\p}\p^{\frac{e_{\p,K}'-1}{2}\frac{n}{e_{\p,K}}}\left(\frac{\prod_{\p}\p^{(e_{\p,L}-1)\frac{n_1(2)}{e_{\p,L}}}}{\alpha}\right)^\frac{n}{2n_1(2)}\]
and we conclude that $x=x_1x_2$, from which we obtain the desired inclusion.
\end{proof}

\begin{teo}\label{abelpari2c}
Let $k$ be a number field, let $G=C(n_1)\times\dots\times C(n_r)$, with $n_{i+1}|n_i$, be an abelian group of order $n$. If $2|n$ and $n_2(2)=1$ then
\[\rt(k,G)=\rt(k,C(n_1(2)))^{\frac{n}{n_1(2)}} \prod_{\substack{l|n\\l\neq 2}}W(k,n_1(l))^{\frac{l-1}{2}\frac{n}{n_1(l)}}.\]
\end{teo}

\begin{proof}
\begin{enumerate}
\item[$\subseteq$] This is Lemma \ref{legameciclichenonr11}.
\item[$\supseteq$] 
This follows by Lemma \ref{multiplicationmetacyclic}, Proposition \ref{konstruktiv} and Lemma \ref{dalleciclichealleabeliane}.
\end{enumerate}
\end{proof}

Thus in any case we reconduct the study of the realizable Steinitz classes for abelian groups to that of $2$-power order cyclic groups. As a consequence of our results we also prove the following corollary.

\begin{coroll}\label{diretto2}
Let $k$ be a number field, let $G$ be an abelian group of order $n$ and let $G_l$ be its $l$-Sylow subgroup for any prime $l|n$. Then
\[\rt(k,G)=\prod_{l|n}\rt(k,G_l)^\frac{n}{n(l)}.\]
\end{coroll}

\begin{proof}
This is immediate by Theorem \ref{oddabelian}, Theorem \ref{abelpari2nc} and Theorem \ref{abelpari2c}.
\end{proof}

In \cite{Cobbe} we prove a similar result concerning a relation between the realizable classes for two groups and for their direct product, in a quite general situation, which however does not include abelian groups of even order.
\end{section} 

\nocite{McCulloh}
\nocite{Long0}
\nocite{Long2}
\nocite{Endo}
\nocite{Carter}
\nocite{Massy}
\nocite{Sodaigui1}
\nocite{Sodaigui2}
\nocite{Soverchia}
\nocite{GodinSodaigui_A4}
\nocite{GodinSodaigui_ottaedri}
\nocite{ByottGreitherSodaigui} 
\nocite{CarterSodaigui_quaternionigeneralizzati}
\nocite{BrucheSodaigui}

\bibliography{bibsteinitz}

\begin{thebibliography}{10}

\bibitem{Bruche}
C.~Bruche.
\newblock Classes de {S}teinitz d'extensions non ab\'eliennes de degr\'e
  {$p^3$}.
\newblock {\em Acta Arith.}, 137(2):177--191, 2009.

\bibitem{BrucheSodaigui}
C.~Bruche and B.~Soda{\"{\i}}gui.
\newblock On realizable {G}alois module classes and {S}teinitz classes of
  nonabelian extensions.
\newblock {\em J. Number Theory}, 128(4):954--978, 2008.

\bibitem{ByottGreitherSodaigui}
N.~P. Byott, C.~Greither, and B.~Soda{\"{\i}}gui.
\newblock Classes r\'ealisables d'extensions non ab\'eliennes.
\newblock {\em J. Reine Angew. Math.}, 601:1--27, 2006.

\bibitem{Carter}
J.~E. Carter.
\newblock Steinitz classes of a nonabelian extension of degree {$p\sp 3$}.
\newblock {\em Colloq. Math.}, 71(2):297--303, 1996.

\bibitem{Carter2}
J.~E. Carter.
\newblock Steinitz classes of nonabelian extensions of degree {$p^3$}.
\newblock {\em Acta Arith.}, 78(3):297--303, 1997.

\bibitem{CarterSodaigui_quaternionigeneralizzati}
J.~E. Carter and B.~Soda{\"{\i}}gui.
\newblock Classes de {S}teinitz d'extensions quaternioniennes
  g\'en\'eralis\'ees de degr\'e {$4p\sp r$}.
\newblock {\em J. Lond. Math. Soc. (2)}, 76(2):331--344, 2007.

\bibitem{tesi}
A.~Cobbe.
\newblock {\em Steinitz classes of tamely ramified Galois extensions of
  algebraic number fields}.
\newblock PhD thesis, Scuola Normale Superiore, Pisa, 2010.

\bibitem{Cobbe}
A.~Cobbe.
\newblock {\em Steinitz classes of tamely ramified Galois extensions of
  algebraic number fields}.
\newblock arXiv:0910.5080v1, to appear in Journal of Number Theory.

\bibitem{Endo}
L.~P. Endo.
\newblock {\em Steinitz classes of tamely ramified Galois extensions of
  algebraic number fields}.
\newblock PhD thesis, University of Illinois at Urbana-Champaign, 1975.

\bibitem{GodinSodaigui_A4}
M.~Godin and B.~Soda{\"{\i}}gui.
\newblock Classes de {S}teinitz d'extensions \`a groupe de {G}alois {$A\sb 4$}.
\newblock {\em J. Th\'eor. Nombres Bordeaux}, 14(1):241--248, 2002.

\bibitem{GodinSodaigui_ottaedri}
M.~Godin and B.~Soda{\"{\i}}gui.
\newblock Module structure of rings of integers in octahedral extensions.
\newblock {\em Acta Arith.}, 109(4):321--327, 2003.

\bibitem{Lang}
S.~Lang.
\newblock {\em Algebraic number theory}.
\newblock GTM 110. Springer-Verlag, New York, second edition, 1994.

\bibitem{Long2}
R.~Long.
\newblock Steinitz classes of cyclic extensions of degree {$l\sp{r}$}.
\newblock {\em Proc. Amer. Math. Soc.}, 49:297--304, 1975.

\bibitem{Long0}
R.~L. Long.
\newblock Steinitz classes of cyclic extensions of prime degree.
\newblock {\em J. Reine Angew. Math.}, 250:87--98, 1971.

\bibitem{Massy}
R.~Massy and B.~Soda{\"{\i}}gui.
\newblock Classes de {S}teinitz et extensions quaternioniennes.
\newblock {\em Proyecciones}, 16(1):1--13, 1997.

\bibitem{McCulloh}
L.~R. McCulloh.
\newblock Cyclic extensions without relative integral bases.
\newblock {\em Proc. Amer. Math. Soc.}, 17:1191--1194, 1966.

\bibitem{McCulloh_Crelle}
L.~R. McCulloh.
\newblock Galois module structure of abelian extensions.
\newblock {\em J. Reine Angew. Math.}, 375/376:259--306, 1987.

\bibitem{Narkiewicz}
W.~Narkiewicz.
\newblock {\em Elementary and analytic theory of algebraic numbers}.
\newblock Springer Monographs in Mathematics. Springer-Verlag, Berlin, third
  edition, 2004.

\bibitem{Neukirch}
J.~Neukirch.
\newblock {\em Class field theory}, volume 280 of {\em Grundlehren der
  Mathematischen Wissenschaften [Fundamental Principles of Mathematical
  Sciences]}.
\newblock Springer-Verlag, Berlin, 1986.

\bibitem{Rotman}
J.~J. Rotman.
\newblock {\em An introduction to the theory of groups}.
\newblock GTM 148. Springer-Verlag, New York, fourth edition, 1995.

\bibitem{Sodaigui1}
B.~Soda{\"{\i}}gui.
\newblock Classes de {S}teinitz d'extensions galoisiennes relatives de degr\'e
  une puissance de 2 et probl\`eme de plongement.
\newblock {\em Illinois J. Math.}, 43(1):47--60, 1999.

\bibitem{Sodaigui2}
B.~Soda{\"{\i}}gui.
\newblock Relative {G}alois module structure and {S}teinitz classes of dihedral
  extensions of degree {$8$}.
\newblock {\em J. Algebra}, 223(1):367--378, 2000.

\bibitem{Soverchia}
E.~Soverchia.
\newblock Steinitz classes of metacyclic extensions.
\newblock {\em J. London Math. Soc. (2)}, 66(1):61--72, 2002.

\end{thebibliography}
\addcontentsline{toc}{section}{Bibliography}
\bibliographystyle{abbrv}

\end{document}